\documentstyle{amsppt}

\nologo
\NoRunningHeads
\TagsOnRight
\NoBlackBoxes

\define\wh{\widehat}
\define\BR{\Bbb R}
\define\BC{\Bbb C}
\define\BT{\Bbb T}
\define\BZ{\Bbb Z}
\define\CF{\Cal F}
\define\CX{\frak X}
\define\CP{\Cal P}
\define\CB{\Cal B}
\define\CL{\Cal L}

\define\CC{\Cal C}
\define\CK{\Cal K}

\define\ve{\varepsilon}
\define\hb{\hbar}
\define\th{\theta}
\define\pa{\partial}
\define\opa{\overline{\partial}}
\define\oz{\overline{z}}

\define\oq{\overline{q}}
\define\oA{\overline{A}}
\define\onu{\overline{\nu}}

\define\od{\overset{\text{def}}\to{=}}
\define\omc{\omega_{\text{\rm class}}}

\define\tg{\operatorname{tg}}
\define\arctg{\operatorname{arctg}}

\define\const{\operatorname{const}}
\redefine\Re{\operatorname{Re}}
\redefine\Im{\operatorname{Im}}

\topmatter
\title Coherent transforms 
and irreducible representations
corresponding to complex structures 
on a cylinder and on a torus
\endtitle

\author M.~V.~Karasev and E.~M.~Novikova\footnotemark\endauthor
\abstract
We study a class of algebras with non-Lie commutation relations
whose symplectic leaves are surfaces of revolution:
a cylinder or a torus.
Over each of such surfaces we introduce 
a family of complex structures and Hilbert spaces of 
antiholomorphic sections 
in which the irreducible Hermitian representations 
of the original algebra are realized.
The reproducing kernels of these spaces 
are expressed in terms of the Riemann theta-function
and its modifications.
They generate quantum K\"ahler structures 
on the surface and 
the corresponding quantum reproducing measures. 
We construct coherent transforms intertwining 
abstract representations of an algebra 
with irreducible representations, 
and these transforms are also expressed 
via the theta-function.
\endabstract

\endtopmatter

\footnotetext{This research was supported by the Russian
Foundation for Basic Research under grant no.~99-01-01047.}

\head 1. Introduction\endhead

A {\it coherent transform} is a linear mapping 
intertwining a given representation of an algebra or 
a group with some irreducible pseudodifferential 
model.
By a pseudodifferential model we mean an algebra 
of pseudodifferential operators 
with symbols on a symplectic manifold. 
We call this manifold a {\it base of a coherent transform}.

The coherent transform is
an important analytical and geometrical tool  
relating objects of classical geometry and quantum objects. 
Moreover, 
this relation is established 
at the level of underlying vector spaces
(not only at the level of the correspondence
{\it symbol\/} $\leftrightarrow$ {\it operator\/}).

The integral kernels of coherent transforms are called 
{\it coherent states}. 
In somewhat different form,
such states appeared already at the dawn
of quantum mechanics 
in the works due to Schr\"odinger,
Heisenberg, Fock,
and then were comprehended by Klauder~\cite{1} 
and Berezin~\cite{2, 3} from the viewpoint 
of the theory of quantization. 
Different generalizations, applications, and references to the
literature devoted to coherent states 
can be found in~\cite{4,~5}.

The above definition of the coherent transform is very general. 
For example, 
it involves constructions of geometric and asymptotic 
quantization \cite{6--8}. 
The construction of coherent transforms 
in such a wide range 
over arbitrary base manifolds is still an open problem.

An important special class of coherent transforms
corresponds to {\it complex structures}
on the base manifold.
In this case the pseudodifferential model 
can be realized 
by using the Wick or Toeplitz operators.
Such operators act in the Hilbert space 
of antiholomorphic sections over the base manifold,
and this space is characterized by its 
{\it reproducing kernel\/} 
(the Bergmann function) \cite{9--11}.
The coherent transform intertwines 
the space of antiholomorphic sections 
with the Hilbert space in which 
the original representation of an algebra is given.

Historically, the first example of coherent transform 
for the three-dimensional Heisenberg algebra 
originates 
precisely from the existence of a complex structure
on the base manifold~$\BR^2$
(here the Gaussian exponential serves as the reproducing kernel,
and the Bargmann coherent transform\cite{12} 
intertwines the Schr\"odinger representation 
and the Fock representation of the Heisenberg algebra).

In such a complex framework, 
the coherent transforms were constructed for 
bounded homogeneous domains in~$\BC^n$ 
and 
for other K\"ahlerial manifolds~\cite{4, 13--15}. 
In the nonhomogeneous case in which 
neither Lie algebras nor Lie groups are present
and representations 
of certain general commutation relations are considered, 
the construction of coherent transforms 
is much less studied, e.g., see~\cite{16--24}.

In the paper~\cite{23}, 
where the approach of \cite{20, 21} was developed,
the authors systematically study 
the non-Lie algebras whose complex coherent transforms 
have, as base manifold, some surfaces of revolution  
(see also generalizations in~\cite{19, 22}). 

It was shown that if the base surface 
(or its closure)
is topologically a plane or a sphere, 
then 
an arbitrary hypergeometric or $q$-hypergeometric function
can appear as the reproducing kernel over this surface.
The correspondence between such special functions, 
various complex structures over the base manifold,
and the representations of a class of algebras 
with non-Lie commutation relations
were established.
With respect to these algebras,
the base manifolds serve as irreducible leaves.

In the present paper 
we continue the study of the algebras from~\cite{23}, 
construct their coherent states over base manifolds
diffeomorphic to the cylinder and the torus, 
and show their relation to the theta-functions.

It should be pointed out that 
the correspondence between complex structures 
on surfaces of revolution 
and special functions, 
which we have constructed,
is 
universal and invariant 
(i.e., is independent of the choice of any 
bases in the spaces considered).
This is a distinction between this correspondence 
and the well-known interpretations of special functions
as matrix elements in the representation theory~\cite{25}
(where the fixing of the basis is important). 
Moreover, we stress that 
the reproducing kernel determines 
the quantum K\"ahler structure over the base manifold
and the quantum reproducing measure,
which in the nonhomogeneous case essentially 
differ from the classical K\"ahler form 
and the classical Liouville measure.
Thus under this approach
the special functions 
(the hypergeometric and theta functions)
generate new geometric objects,
which play the key role in quantum geometry 
of surfaces of revolution.

Needless to say that,
such a base manifold as a torus
possessing a pair of noncontractible $1$-cycles
is always naturally associated with noncommutative discrete groups.
In particular,  
the base manifolds is associated with discrete 
subgroups of the Heisenberg group, 
which, as is well-known, 
are closely related to the theory 
of theta-functions~\cite{26--28}.  
This fact has been used in a number of papers 
including those where the completeness 
of coherent states over elementary cells of the plane 
was studied~\cite{4} 
and those related to the Weyl quantization over the
torus~\cite{29}.  
However, our results are significantly different.
In particular, 
we associate the irreducible representations 
of a noncommutative algebra and the Riemann theta-function,
first of all, with a cylinder on which there is only one cycle.
In the case of a torus, 
our correspondence implies already the product of two
theta-functions of two different arguments\footnotemark. 

\footnotetext{The latter correlates with Sklyanin's
hypothesis~\cite{16}, which, though,  
deals with a somewhat different situation.}

Finally, we mention one more (different) correspondence 
between special functions and noncommutative algebras
which accompanies 
the Yang--Baxter and Kni\-zhnik--Zamolodchikov equations
and is related to the parametrization 
of the set of commutation relations 
(structure constants), e.g., see~\cite{16,~30}.

The results of the present paper
were partially presented in the report 
of one of the authors 
at the Euroconference dedicated to the memory of M.~Flato
(Dijon, September 2000); see~\cite{31} 
where a general discussion of quantum geometry 
and certain applications are given. 
The first of the authors is very grateful 
to A.~Weinstein, P.~Cartier, and  A.~A.~Kirillov
for useful remarks.

\head 2. Commutation relations and surfaces of revolution\endhead

We consider a noncommutative algebra with involution 
and with Hermitian generators 
$\wh{S}_1,\wh{S}_2,\wh{A}_1,\dots,\wh{A}_k$
satisfying several commutation relations.
To describe these relations,
we introduce the following complex combinations of generators: 
$$
\wh{B}=\wh{S}_1-i\wh{S}_2,\qquad 
\wh{C}=\wh{S}_1+i\wh{S}_2
$$
and choose a one-parametric group of transformations
$$
\Phi_t\:\BR^{k+1}\to\BR^{k+1},\qquad -\infty<t<\infty
$$
(i.e., the flow of a vector field).
For the coordinates in $\BR^{k+1}$ we write $A_0,A_1,\dots,A_k$
and denote the last~$k$ coordinates by $A=(A_1,\dots,A_k)$. 
Then $\Phi_t$ has the form
$$
\Phi_t(A_0,A)=\big(\varphi^0_t(A_0,A),\varphi_t(A_0,A)\big).
\tag 2.1
$$
Here $\varphi^0_t$ is a scalar function 
(the zero component of $\Phi_t$),
and $\varphi_t$ is a vector function ranging in~$\BR^k$. 
In this notation we have 
$$
\varphi^0_t(A_0,A)\big|_{t=0}\equiv A_0,\qquad 
\varphi_t(A_0,A)\big|_{t=0}\equiv A.
\tag 2.2
$$

Into all components of the function~(2.1), 
instead of the coordinates $A_0$, $A$,
we can substitute Hermitian elements of the algebra.
We postulate the following relations between the generators 
of the algebra:
$$
\gathered
\wh{C}\cdot\wh{B}=\varphi^0_\hb(\wh{B}\wh{C},\wh{A}),\qquad 
\wh{C}\cdot\wh{A}=\varphi_\hb(\wh{B}\wh{C},\wh{A})\cdot\wh{C},
\qquad \wh{A}_j\cdot \wh{A}_l=\wh{A}_l\cdot\wh{A}_j,\\
\wh{B}^*=\wh{C}, \qquad \wh{A}^*_j=\wh{A}_j,
\endgathered
\tag 2.3
$$
where $\hb>0$ is a fixed number and $j,l=1,\dots,k$. 
After conjugation, 
the second relation in (2.3) also implies 
$$
\wh{A}\cdot\wh{B}=\wh{B}\cdot \varphi_\hb(\wh{B}\wh{C},\wh{A}).
\tag 2.4
$$

Note that, by (2.3) and (2.4), 
the element $\wh{B}\wh{C}$ commutes with all elements 
$\wh{A}=(\wh{A}_1,\dots,\wh{A}_k)$ and hence the functions 
of $\wh{B}\wh{C}$ and $\wh{A}$ 
on the right-hand sides in~(2.3) and~(2.4) are well defined.

Formulas (2.2) show that the algebra with relations~(2.3)
is a deformation with respect to the parameter~$\hb$ 
of a commutative algebra (corresponding to $\hb=0$).
Note that the parameter $\hb$ replaces 
the ``time''-argument~$t$ in~(2.3) and~(2.4).
Hence our noncommutative algebra is generated 
by the ``deformation'' flow~$\Phi_t$. 
A similar idea to take~$\hb$ as time 
was used in~\cite{29} to derive the differential equation 
for the Weyl $*$-product, although the deformation flow 
was not considered there.

Obviously, the standard {\it Casimir elements\/} 
of the algebra with relations~(2.3)
(which enter the center for any $\hb>0$) 
have the form\footnotemark
$$
\varkappa(\wh{B}\wh{C},\wh{A}),
\tag 2.5
$$
where the function $\varkappa$ is constant on the trajectories 
of the flow~$\Phi_t$,
i.e., $\varkappa(\Phi_t)=\varkappa$. 
The number of such independent functions $\varkappa$
is equal to~$k$.
We are going to study 
{\it operator irreducible representations\/}
of the algebra~(2.3). 
In each such representation all Casimir elements 
are scalar, and we obtain the following~$k$ relations:
$$
\varkappa_j(\wh{B}\wh{C},\wh{A})=\const_j\cdot I,\qquad j=0,1,\dots,k-1.
\tag 2.6
$$
Instead of these equations, 
we can introduce the elements $\wh{B}\wh{C}$ and $\wh{A}$
as functions of a single Hermitian element $\wh{t}=\wh{t}^{*}$:
$$
\wh{B}\wh{C}=\varphi^0_{\wh{t}}(a_0,a),\qquad
\wh{A}=\varphi_{\wh{t}}\,(a_0,a),
\tag 2.7
$$
where $(a_0,a)$ is a chosen point in  $\BR^{k+1}$
such that $a_0\geq 0$ and
$\varkappa_j(a_0,a)=\const_j$, $j=0,\dots,k-1$. 

\footnotetext{In general, 
the center of the algebra~(2.3) 
is not exhausted by such elements; 
in what follows, see comments 
on the quantization condition~(3.12) 
and Example~5.1.}

By using Eqs.~(2.7) or (2.6),
we assign to the quantum relations 
the classical surface $\CX\subset \BR^{k+2}$:
$$
\aligned
\CX&=\{BC=\varphi^0_t(a_0,a),\; A=\varphi_t(a_0,a)\}\qquad
\text{or}\\
\CX&=\{\varkappa_j(BC,A)=\varkappa_j(a_0,a),\; j=0,\dots,k-1\},
\endaligned
\tag 2.8
$$
where $B\equiv\overline{C}\equiv S_1-iS_2$.
Since $BC=S^2_1+S^2_2$, 
for each chosen~$t$
Eqs.~(2.8) determine 
a circle of radius $\varphi^0_t(a_0,a)^{1/2}$
in the plane of variables $S_1,S_2$.
If $t$ varies, these circles fill a surface of revolution 
embedded in~$\BR^{k+2}$. 

The geometry of the surface depends on the properties 
of the function $\varphi^0_t(a_0,a)$. 
For $a_0=0$, 
if $\varphi^0_t(0,a)>0$ either for all $t>0$ 
or for all $t<0$, 
then the closure of the surface $\CX$ is topologically 
equivalent to the plane.
However, if either $\varphi^0_t(0,a)>0$ only for $0<t<t_0$
and $\varphi^0_{t_0}(0,a)=0$
or  $\varphi^0_t(0,a)>0$ for $-t_0<t<0$ 
and $\varphi^0_{-t_0}(0,a)=0$, 
then the closure of the surface $\CX$ 
is topologically equivalent to the sphere.
Both these cases are described in~\cite{23}.

In the present paper we are interested in the case $a_0>0$
and $\varphi^0_t(a_0,a)>0$ for all $t\in\BR$. 
Then~$\CX$ is either a cylinder $\Bbb S\times\BR$
embedded in~$\BR^{k+2}$
or if, in addition, we have the periodicity  
$$
\exists T:\quad
\varphi^0_{_{\ssize T}}(a_0,a)=a_0,\qquad 
\varphi^{}_{_{\ssize T}}(a_0,a)=a,
\tag 2.9
$$
then the surface~$\CX$ is a torus $\BT^2$
embedded in~$\BR^{k+2}$. 

We construct irreducible representations and coherent states 
of the algebra~(2.3) which correspond to complex structures 
on the cylinder or on the torus $\CX\subset\BR^{k+2}$.

\head 3. Representations and complex structures\endhead
  
We introduce another Hermitian generator 
$\wh{s}={\wh{s}\,}^{*}$ 
whose canonical commutation relation with~$\wh{t}$
is $[\,\wh{t},\wh{s}\,]=-i\hb\cdot I$
and which is related to~$\wh{B}$ and $\wh{C}$ 
by the formulas
$$
\wh{B}=\mu(\wh{t}\,)\exp\{i\wh{s}\},\qquad
\wh{C}=\exp\{-i\wh{s}\}\overline{\mu}(\wh{t}\,).
\tag 3.1
$$
Here the complex function~$\mu$ 
(determining the modulus of~$B$ and~$C$
as well as the origin point of the argument~$s$) 
is subject to the condition
$$
|\mu(t)|^2=\CF(t),\qquad\text{where}\quad
\CF(t)\od\varphi^0_t(a_0,a)>0.
\tag 3.2
$$ 

Thus formulas~(3.1) introduce 
Darboux coordinates on the surface $\CX$~(2.8).
In addition, the surface is $2\pi$-periodic 
with respect to the coordinate~$s$,
i.e.,  $s$ corresponds to a noncontractible cycle on~$\CX$. 
We want to introduce a complex structure on~$\CX$
by using the Darboux coordinates.

Assume that we realized the operator~$\wh{t}$ 
as the derivation by a complex coordinate, 
i.e., 
$\wh{t}=\hb\pa/\pa \oz$, $\oz\in\BC$.
Then, in the space of functions of~$\oz$,
relations~(2.3) and~(2.4)
can be realized formally by the operators 
$$
\wh{B}=\CB(\wh{t}\,)\cdot\exp\{-\tau\wh{t}+\wh{\oz}\},\qquad 
\wh{C}=\exp\{\tau\wh{t}-\wh{\oz}\}\cdot\CC(\wh{t}\,),\qquad
\wh{A}=\varphi_{\wh{t}}\,(a_0,a),
\tag 3.3
$$
where $\wh{\oz}$ is the multiplication operator by~$\oz$,
$\tau>0$ is a constant, and the functions~$\CB$ and~$\CC$ 
satisfy the condition
$$
\CF(t)=\CB(t)\CC(t).
\tag 3.4
$$

Relations (3.3) and (3.1) easily imply
$$
\exp\{\wh{\oz}\}=\exp\{\tau\wh{t}+g(\,\wh{t}\,)+i\wh{s}\}.
\tag 3.5
$$
Here $g$ is a solution\footnotemark 
of the equation
$$
\exp\bigg\{\frac1{\hb}\int^{t}_{t-\hb}g(t)\,dt\bigg\}
=\nu(t),\qquad \text{where}\quad
\nu \od \frac{\mu}{\CB}.
\tag 3.6
$$
The quantum relation (3.5) prompts us to assign 
the complex coordinate~$\oz$ to the Darboux coordinates~$t$ and~$s$ 
as follows:
$$
\oz=\tau t+g(t)+is.
\tag 3.7
$$
The mapping defined by this formula is one-to-one if 
$$
\tau>\tau_0,\qquad 
\tau_0\od -\min_t\Re\big(g'(t)\big)<\infty.
\tag 3.8
$$
It is easy to calculate that 
$$
\aligned
&\frac{\pa}{\pa\oz}=\frac12\bigg(\frac{1}{\tau+\Re g'(t)}\,
\frac{\pa}{\pa t}-i\frac{\tau+\overline{g'(t)}}{\tau+\Re g'(t)}\, 
\frac{\pa}{\pa s}\bigg),\\
&\omc\od dt\wedge ds
=\frac{i}{2(\tau+\Re g')}d\oz\wedge dz.
\endaligned
\tag 3.9
$$

\footnotetext{This is the solution that is regular as $\hb\to 0$; 
the general solution of~(3.6) allows us
to add an arbitrary $\hb$-periodic function with zero mean value.}

\proclaim{Lemma 3.1}
Suppose that the condition {\rm(3.8)} is satisfied
and in the periodic case {\rm(2.9)}, in addition, 
the function~$g'$ is $T$-periodic.
Then formulas {\rm(3.9)} define 
a complex structure and a K\"ahler form on~$\CX$.
\endproclaim

First we study the case in which the surface~$\CX$ 
is homeomorphic to the cylinder.
The estimate $\tau_0<\infty$ (3.8) can be provided as follows.
By~$S_m$ we denote the class of smooth functions~$g$ 
on the straight line 
whose derivatives $g^{(k)}$ satisfy the estimate
$$
\exists r\le m \quad\forall k\ge1\quad \exists c\:
|g^{(k)}(t)|\leq c(1+|t|)^{r-k}.
$$

\proclaim{Lemma 3.2}
Let the function $\nu=\mu/\CB$ be such that $\nu'/\nu\in S_{-1}$.
Then there exists a unique solution $g\in S_1$ of~{\rm(3.6)}. 
Since $g'(t)$ is bounded for $|t|\to\infty$,
the estimate~{\rm(3.8)} holds for some~$\tau_0$.
\endproclaim

Now we write explicit formulas for the solution of~(3.6).
Assume that $\nu$ has the asymptotics
$$
\nu(t)=|t|^b e^{p t+l }\big(1+O(t^{-1})\big)
$$
either as $t\to-\infty$ or as $t\to+\infty$. 
Then the solution $g\in S_1$ of~(3.6) has the form
$$
g(t)=\sum^{\infty}_{k=1}\Big(\hb\frac{\nu'}{\nu}(t-k\hb)
+\frac{b}k-\hb p\Big)
+\hb\frac{\nu'(t)}{\nu(t)}+l+b(\ln\hb-\gamma)
+p\Big(t-\frac{\hb}{2}\Big)
$$
or respectively
$$
g(t)=\sum^{\infty}_{k=1}\Big(-\hb\frac{\nu'}{\nu}(t+k\hb)
+\frac{b}k+\hb p\Big)+l+b(\ln\hb-\gamma)
+p\Big(t+\frac{\hb}{2}\Big).
$$
Here $\gamma$ stands for the Euler constant
$\gamma=\lim_{m\to\infty}\big(1+\frac12+\dots+\frac1m-\ln m\big)=
0,577\dots\,\,$.

\example{Remark 3.1}
Under the special choice of the factor $\CB=\mu$
we obtain $\nu=1$ and $g=0$. 
The following special complex structure corresponds to this
case: 
$$
\oz=\tau t+is,\qquad 
\omc=\frac{i}{2\tau}d\oz\wedge dz.
\tag 3.7a
$$
Conversely, 
an arbitrary complex structure of the form (3.9) 
can be obtained by formulas (3.6) 
by choosing an appropriate factor~$\CB$.  
Thus the family of complex structures on~$\CX$ 
is always determined by the factorization choice~(3.4) 
for the zero component $\CF(t)=\varphi^0_t(a_0,a)$ 
of the trajectory of the flow~$\Phi_t$.
\endexample

Now we prove that the formal operator realization~(3.3) 
becomes 
the actual irreducible representation of relations~(2.3).

For any nonzero function $\nu$ on the real line 
we define
$$
\nu_!(n\hb)=\cases
\nu(\hb)\cdot\nu(2\hb)\cdot\dots\cdot\nu(n\hb),\quad&n\ge1,\\
1,\quad&n=0,\\
\nu(0)^{-1}\cdot\nu(-\hb)^{-1}\cdot\dots\cdot\nu((n+1)\hb)^{-1},
\quad&n\le-1.
\endcases
$$
Now let us consider the function~$\nu$ defined by~(3.6), 
i.e.,
$\nu=\mu/\CB$
and introduce the Hilbert space $\CL_{\nu}$
of $2\pi i$-periodic antiholomorphic functions 
on the complex plane 
equipped with the norm
$$
\|\psi\|=\Big(\sum_{n\in \BZ}
|\nu_{!}(n\hb)|^2 e^{\tau\hb n^2}|\psi_n|^2\Big)^{1/2}\quad  
\text{if}\quad
\psi(\oz)=\sum_{n\in\BZ}\psi_n e^{n\oz}.
\tag 3.10
$$

Under the conditions of Lemma~3.2,
if the inequality $\tau>\tau_0$~(3.8) holds, 
then the space $\CL_\nu$ is identified 
with the space of antiholomorphic functions 
on the cylinder~$\CX$ with complex structure~(3.9).

\proclaim{Theorem 3.1}
Suppose that the periodicity condition~{\rm(2.9)} 
does not hold \rom(the cylinder case\rom).
Then formulas {\rm(3.3)} 
determine an infinite-dimensional, 
irreducible in the operator sense,
representation of relations~{\rm(2.3)} 
in the Hilbert space~$\CL_\nu$ 
of antiholomorphic functions on the cylinder.
A different choice of the factors~$\CB$ and~$\CC$ 
\rom(and hence the choice of the complex structure 
on the cylinder\/\rom), 
as well as a different choice of the parameter~$\tau$
leads to equivalent representations.
\endproclaim

Now we study the case of a torus.
Suppose that the $T$-periodicity condition~(2.9)
is satisfied. 
It should be noted that the quantum generator~$\wh{t}$
in the irreducible representation has a discrete spectrum
consisting of the numbers~$n\hb$, where $n$ is integer. 
Hence, there are two versions in the periodic case:
either the period~$T$ with respect to the classical
coordinate~$t$ is not a multiple of~$\hb$, 
i.e.,
$$
T\not=N\hb\qquad \text{not for any integer~$N$},
\tag 3.11
$$ 
or $T$ is a multiple of~$\hb$ 
and then the condition $T=N\hb$
can be written as an integral over the torus
(or over its finite sheet covering):
$$
\frac1{2\pi \hb}\int_{\CX}\omc=N.
\tag 3.12
$$
Let us explain this in more detail.
It is not assumed that the period~$T$ in~(2.9) is minimal.
Let $T_0$ be the minimal period.
Then the nonresonance version (3.11) means that~$T_0$
and~$\hb$ are incommensurable.
Conversely, 
in the resonance version 
we have  $T_0=N\hb/m$, where~$N$ and~$m$ are coprime 
integers. 
Then the $m$-multiple period $T=mT_0$ 
satisfies the condition $T=N\hb$
and thus the integral in~(3.12) 
must be taken over the $m$-sheet covering of the torus~$\CX$. 

{\it In the nonresonance version {\rm(3.11)} 
Theorem~{\rm 3.1} remains valid},
i.e, the representation~(3.3) is irreducible in the operator
sense. 
We point out that this representation is infinite-dimensional, 
although it corresponds to a compact manifold! 
The cylinder in Theorem~3.1  
is the infinite sheet covering of the torus 
(the infinite winding proceeds along the $t$-axis).

In the resonance version (3.12),
the representation~(3.3) in the space~$\CL_\nu$ 
is vector irreducible but already not operator irreducible.
Indeed, 
in this case~$\wh{B}^N$ is not a scalar operator 
but it commutes with all operators of the representation.
In the  representations theory, 
such 
``nonclassical'' Casimir operators are not quite usual objects. 
For instance, 
in the paper~\cite{16} where 
the irreducible representations corresponding to tori
were studied for the first time, 
these Casimir operators were not presented.
We point out that, 
in the classical limit as $\hb\to0$, 
there does not exist any function 
on the Poisson manifold 
that corresponds to such elements, 
since $B^{N}=B^{T/\hb}$ does not have any limit 
as $\hb\to0$. 
Later on
this effect is considered in 
Example~5.1. 

\example{Remark 3.2}
The complex structure on the torus $\CX$ 
is introduced by the general formulas~(3.9), 
where $g'$ is now a $T$-periodic function.
We also assume that the function~$\mu$ in~(3.2) 
is $T$-periodic and the factors~$\CB$ and~$\CC$ in~(3.4) 
are determined by~$g$ and~$\mu$
by the formulas~(3.6).  

Here the inverse problem of constructing~$g$ 
either from a given factor $\CB$ 
or from the function $\nu=\mu/\CB$ 
(the problem we considered in the nonperiodic case) 
is made more complicated 
in view of the problem of small 
denominators in~(3.6) near the resonance $T\approx N\hb$.
At the resonance $T=N\hb$
we have  
additional conditions that the right-hand side 
is orthogonal to some Fourier harmonics.
\endexample

Now let us consider the {\it resonance version}~(3.12) 
and show how 
the operator irreducible representations 
of relations~(2.3) can be constructed in this case. 
First we note that, 
by choosing appropriate constants 
to normalize the $N\hbar$-periodic functions~$\CB$ and~$\CC$
in~(3.4),   
it is possible to make relation~(3.4) still hold and,
simultaneously, to have 
$$
|\nu_{!}(N\hb)|=1,\qquad 
\sum^{N}_{n=1}\arg\CB(n\hb)=\alpha \quad (\bmod 2\pi),
\tag 3.13
$$
where $\alpha$ is a given real number.

Then we introduce a Hilbert space of functions over the
torus~$\CX$. The complex coordinate~$z$ 
identifies the torus with the rectangle in the complex plane
$$
0\leq\Im z<2\pi,\qquad
0\leq\Re z<\tau T=\tau\hb N.
$$
Since there do not exist 
double-periodic antiholomorphic functions 
with periods~$2\pi i$ and~$\tau \hb N$ 
(except constants), 
we use quasiperiodic functions.
We also consider the conditions of quasiperiodicity
$$
\psi(\oz+2\pi i)=\psi(\oz),\qquad
\psi(\oz+\tau\hb N)=
\exp\Big\{\frac{\tau\hb}{2}N^2+N\oz\Big\}\psi(\oz).
\tag 3.14
$$
These conditions can be rewritten as 
$$
\exp\Big\{\frac{2\pi i}{\hb}\wh{t}\Big\}\psi(\oz)
=\psi(\oz),
\qquad 
\exp\Big\{\frac{T}{\hb}(\wh{\oz}-\tau\wh{t}\,)\Big\}\psi(\oz)
=\psi(\oz).
$$
Thus 
the quasiperiodicity is 
the periodicity of powers of the exponents 
of the creation operator $\wh{\oz}-\tau\wh{t}$
and the operator $\frac{2\pi i}{\hb}\wh{t}$.

The functions satisfying the quasiperiodicity conditions 
are the theta-functions.
Let us choose one of them, e.g., 
$$
\th (\alpha,\ve)=\sum_{n\in\BZ}
\exp\{-\ve n^2+in\alpha\},\qquad \ve>0.
\tag 3.15
$$
Then system (3.14) has the solution 
$\psi(\oz)=\th (N\oz/i,\tau\hb N^2/2)$. 
The other solutions can be obtained by applying 
the powers of the exponent of the creation operator 
$\exp\big\{n(\wh{\oz}-\tau\wh{t}\,)\big\}$, $n=0,1,\dots,N-1$.

Thus we arrive at the following description 
of antiholomorphic functions satisfying~(3.14):
$$
\psi(\oz)=\sum^{N-1}_{n=0}\psi_n e^{n\oz}
\th \Big(\frac{N(\oz-\tau\hb n)}{i},\frac{\tau\hb N^2}{2}\Big).
\tag 3.16
$$
We define their norm as 
$$
\|\psi\|=\Big(\sum^{N-1}_{n=0}|\nu_{!}(n\hb)|^2 
e^{\tau\hb n^2}|\psi_n|^2\Big)^{1/2}.
\tag 3.17
$$
The Hilbert space thus obtained
we denote by $\CL^{N}_{\nu}$. 
Under the condition (3.12),
the space $\CL^{N}_{\nu}$ is identified with the space 
of antiholomorphic sections over the torus~$\CX$.

\proclaim{Theorem 3.2}
Suppose that the periodicity condition {\rm(2.9)} 
holds for $T=N\hb$ and not for $T=N'\hb$ 
for any integer $0<N'<N$.
Suppose also that the $T$-periodic factors~$\CB$
and~$\CC$ in~{\rm(3.4)} satisfy~{\rm(3.13)}.
Then formulas~{\rm(3.3)} determine 
an $N$-dimensional operator irreducible representation
of relations~{\rm(2.3)} 
in the Hilbert space~$\CL^{N}_{\nu}$ 
of antiholomorphic sections over the torus
and 
$$
\wh{B}^N=\CF_{!}(N\hb)^{1/2}e^{i\alpha}\cdot I.
\tag 3.18
$$
The numbers $N$ and $\alpha$ are parameters 
of the representation, i.e.,
to different pairs $(N,\alpha)$ 
there correspond nonequivalent representations.
A different choice of the factors~$\CB$ and~$\CC$ 
\rom(and thus a choice of the complex structure on the torus\rom)
as well as a different choice of the parameter~$\tau$
leads to equivalent representations.
\endproclaim

It is natural to pose the question:
does this construction include 
all irreducible representations of relations~(2.3)?

\proclaim{Theorem 3.3}
All operator irreducible representations of~{\rm(2.3)}
for which the operators~$\wh{B}\wh{C}$ 
and $\wh{A}_1,\dots,\wh{A}_k$ have a nonempty point spectrum
can be classified as\/{\rm:}

$(1^0)$ either the eigenvalues $\wh{B}\wh{C}$ and $\wh{C}\wh{B}$ 
are positive and then the representation is equivalent 
to one of those in Theorem~{\rm 3.1} or~{\rm 3.2;} 
this type of representation corresponds 
to the cylinder or the torus~$\CX${\rm ;}

$(2^0)$ or the operator~$\wh{B}\wh{C}$ or $\wh{C}\wh{B}$ 
possess the zero eigenvalues and then the representation 
is equivalent to one of those given in~\cite{23}{\rm;} 
this type of representation corresponds 
to the plane or the sphere~$\CX$.
\endproclaim

\head 4. Reproducing kernels and coherent states\endhead

The Hilbert space of antiholomorphic functions can be
characterized by its reproducing kernel.
First, we study the case of a cylinder.
The reproducing kernel corresponding to the space~$\CL_\nu$
has the form 
$\CK_\nu=\sum_n|e^{(n)}|^2$,
where $\{e^{(n)}\}$ is an orthonormal basis in~$\CL_\nu$.
Starting from~(3.10), 
we choose the basis
$$
e^{(n)}(\oz)=\nu_{!}(n\hb)^{-1}
\exp\Big\{-\frac{\tau\hb n^2}{2}+n\oz\Big\},\qquad 
n\in\BZ.
\tag 4.1
$$
Then for $\CK_\nu$ we obtain the series
$$
\CK_\nu(\oz|z)
=\sum_{n\in\BZ}\frac1{|\nu_{!}(n\hb)|^2}
\exp\{-\tau\hb n^2+n(\oz+z)\}
=\th_{|\nu|^2}\Big(\frac{\oz+z}{i},\tau\hb\Big).
\tag 4.2
$$
Here by~$\th_\rho$ stands for the following modification 
of the theta-function:
$$
\th_\rho(\alpha,\ve)\od\sum_{n\in\BZ}
\frac1{\rho_{!}(n\hb)}\exp\{-\ve n^2+in\alpha\}
=\rho_{!}\Big(-i\hb\frac{d}{d\alpha}\Big)^{-1}
\th (\alpha,\ve),
$$
where $\th$ is the standard theta-function (3.15). 
Certainly, here and in~(4.2)
some estimates for $\rho$ or $|\nu|$ at $\pm\infty$,
which ensure the uniform convergence of infinite series,
are assumed to be satisfied. 

Note that the modified theta-function~$\th_\rho$
in terms of which the reproducing kernel is expressed 
can be defined as the unique $2\pi$-periodic solution 
of the problem 
$$
\exp\Big\{i\Big(2\ve\frac{d}{d\alpha}+\alpha\Big)\Big\}y(\alpha)
=\rho\Big(-i\hb\frac{d}{d\alpha}\Big)y(\alpha),
\qquad
\frac1{2\pi}\int^{2\pi}_{0}y(\alpha)\,d\alpha=1.
\tag 4.3
$$

Next, the imaginary Jacobi transformation is known for~$\th$:
$$
\th (\alpha,\ve)=\sqrt{\frac{\pi}{\ve}}
\exp\Big\{-\frac{\alpha^2}{4\ve}\Big\}
\th \Big(\frac{i\pi\alpha}{\ve},\frac{\pi^2}{\ve}\Big),
$$
which generates some transformation of the function~$\th_\rho$.
Hence for the reproducing kernel~(4.2) 
we obtain the identity
$$
\CK_\nu(\oz|z)=\sqrt{\frac{\pi}{\tau\hb}}
\exp\Big\{\frac{(\oz+z)^2}{4\tau\hb}\Big\}
q_\nu(\oz+z),
$$
where
$$
q_\nu(x)\od\Big|\nu_{!}\Big(\frac{x}{2\tau}+\hb\frac{d}{dx}\Big)\Big|^{-2}
\th \Big(\frac{\pi x}{\tau\hb},\frac{\pi^2}{\tau\hb}\Big).
$$
Here we use the function $\nu_{!}(t)$ which is defined not only 
at the lattice points $t=n\hb$ but also for arbitrary complex~$t$ 
by the difference equation 
$$
\nu_{!}(t+\hb)=\nu(t+\hb)\cdot \nu_{!}(t),\qquad
\nu_{!}(0)=1
$$
(an analog of the gamma-function). 

Note that the function $\nu_{!}$ depends on the
deformation parameter~$\hb$ in a singular way  
and can be expressed in terms of the solution~$g$
of~(3.6) as 
$$
\align
\nu_{!}(t)
&=\exp\bigg\{\frac{1}{\hb}\int^t_0 g(\tau)\,d\tau\bigg\}\\
&=\sqrt{\frac{\nu(t)}{\nu(0)}}\exp\bigg\{\frac1{\hb}
\int^t_0\ln\nu(\tau)\,d\tau\bigg\}
\prod^{\infty}_{k=1}\exp\Big\{\hb^kb_{k+1}
\big[(\ln\nu)^{(k)}(t)-(\ln\nu)^{(k)}(0)\big]\Big\},
\endalign
$$
where $b_k$ are the Bernoulli numbers:
$x/(1-e^{-x})=\sum^{\infty}_{k=0}x^kb_k$. 
More details formulas for $\nu_!(t)$ 
can be derived by using the series for~$g(t)$
given after Lemma~3.2.

Now we construct one more function
$$
p_\nu(x)\od\Big|\nu_{!}\Big(\frac{x}{2\tau}-\hb\frac{d}{dx}\Big)\Big|^21
=\frac1{\sqrt{\pi\hb\tau}}
\int^{\infty}_{-\infty}\exp\Big\{-\frac{t^2}{\hb\tau}\Big\}
|\nu|_{!}\Big(\frac{x+2it}{2\tau}\Big)^2\,dt.
\tag 4.4
$$
Here it is assumed that the integral in~(4.4) 
converges.
Unfortunately, 
at present we cannot present general sufficient conditions 
on $\nu(t)$ which guarantee the convergence.
However, in examples the convergence is justified 
(e.g., see Example~5.2 below). 

By using the functions $q_\nu$ and $p_\nu$,
we specify the following measure on the cylinder~$\CX$:  
$$
dm_\nu\od(q_\nu p_\nu)(\oz+z)\frac{d\oz\,dz}{2\tau}.
$$

\proclaim{Lemma 4.1}
If the function $p_\nu$ satisfies the estimate 
$|p_\nu(x)|\leq c\exp\{x^2/4\hb\tau\}(1+|x|)^{-1-\ve}$, 
where $\ve>0$, 
then the norm in the space $\CL_\nu$~{\rm(3.10)} 
can be written in the integral form 
$$
\align
\|\psi\|^2&=\frac1{2\pi}\int_{0\leq\Im z\leq 2\pi}
|\psi(\oz)|^2 p_\nu(\oz+z)
\frac{\exp\{-(\oz+z)^2/(4\hb\tau)\}}{\sqrt{4\pi\hb\tau}}\,d\oz\,dz\\
&
=\frac1{2\pi\hb}\int_{\CX}\frac{|\psi|^2}{\CK_\nu}\,dm_\nu.
\tag 4.5
\endalign
$$
\endproclaim

The last expression in (4.5) is the general geometric form 
of the norm in the space of antiholomorphic sections 
over the K\"ahlerian manifold~$\CX$ with K\"ahler form
$$
\omega_\nu\od
i\hb\opa\pa(\ln \CK_\nu)\,d\oz\wedge dz.
\tag 4.6
$$
This form is called a {\it quantum K\"ahler form}, 
and the measure~$dm_\nu$ in (4.5) is called 
a {\it reproducing measure}..

In particular, in the special case (3.7a), i.e., 
for $\nu=1$ we obtain 
$q_1(x)=\th (\pi x/\tau \hb, \pi^2/\tau\hb)$
and $p_1(x)=1$.
In this case the quantum K\"ahler form $\omega=\omega_1$ 
and the reproducing measure $dm=dm_1$ 
on the cylinder $\CX$ 
are determined by the standard theta-function as 
$$
\omega=i\bigg(\frac1{2\tau}+\hb\opa\pa\ln\th 
\Big(\frac{\pi(\oz+z)}{\tau\hb},\frac{\pi^2}{\tau\hb}\Big)\bigg)
d\oz\wedge dz,\quad
dm=\th \bigg(\frac{\pi(\oz+z)}{\tau\hb},\frac{\pi^2}{\tau\hb}\bigg)
\frac{d\oz\,dz}{2\tau}.
\tag 4.7
$$
Of course, in the limit $\hb=0$ 
these quantum geometric objects 
become the classical form and the classical measure
on the cylinder:
$$
\omc=\frac{i}{2\tau}d\oz\wedge dz, \qquad
dm_{\text{\rm class}}=\frac1{2\tau}d\oz\,dz.
$$
However, as $\hb\ne 0$ 
the quantum K\"ahler form and the quantum measure
significantly differ form the classical ones.
Moreover, the quantum reproducing measure does not
coincide with the Liouville measure 
generated by the quantum K\"ahler form.
It is of interest that 
this distinction is exponentially small:
$$
\omega=\omc+O\big(e^{-\pi^2/\tau\hb}\big), \quad
dm=dm_{\text{\rm class}}+O\big(e^{-\pi^2/\tau\hb}\big),\qquad
\hb\to 0.
$$

In conclusion, 
we write the formula for coherent states.
The {\it fiducial state\/} is defined to be
a vector~$\CP^0$ in the space of the representation 
on which the operators $\wh{A}$ and $\wh{A}^0=\wh{B}\wh{C}$ 
take the values~$a$ and~$a_0$ 
(which corresponds to the initial point $t=0$), 
i.e., 
$$
\wh{A}\CP^0=a\CP^0,\qquad 
\wh{B}\wh{C}\,\CP^0=a_0\CP^0,\qquad 
\|\CP^0\|=1.
\tag 4.8
$$
Note that 
the fiducial state of the representation (3.3) 
in the Hilbert space~$\CL_\nu$
is the unit function 
$\CP^0(\oz)\equiv 1\in\CL_\nu$.
Then, by (3.5) and (3.6), we have
$$
e^{n\oz}=e^{n\wh{\oz}} 1=
\nu_{!}(n\hb)\exp\{n^2\tau\hb/2+in\wh{s}\}\CP^0.
$$
Substituting this expression into the formula~(4.2) 
for $\CK_\nu$, we obtain
$$
\CK_\nu=\sum_n\onu_{!}(n\hb)^{-1}
\exp\Big\{-\frac{\tau\hb n^2}{2}+in\wh{s}+nz\Big\}\CP^0
=\th_{\onu}\Big(\frac{z}{i}+\wh{s},\frac{\tau\hb}{2}\Big)\CP^0.
\tag 4.9
$$
This formula determines the coherent states in the
space~$\CL_\nu$. 
If now 
we consider an abstract Hilbert space~$L$ 
where some representation of the algebra~(2.3) acts,
determine the fiducial state $\CP^0\in L$ 
by formulas (4.8), 
and introduce the operator~$\wh{s}$ in~$L$ 
by formulas~(3.1), 
then 
the {\it coherent states\/} $\CP_z\in L$ 
are determined by the same formula~(4.9): 
$$
\CP_z=\th_{\onu}\Big(\frac{z}{i}+\wh{s},
\frac{\tau\hb}{2}\Big)\CP^0.
\tag 4.10
$$

By~$\Pi(\oz|z)$ we denote 
the orthogonal projector in $L$ 
on the subspace generated by~$\CP_z$.

\proclaim{Theorem 4.1}

{\rm(a)} The inner product of coherent states 
{\rm(4.10)} over the cylinder~$\CX$ 
yields the reproducing kernel{\rm:}
$\CK_\nu(\oz|z)=(\CP_z,\CP_z)$.

{\rm(b)} The partition of unity takes place{\rm:}
$$
\frac1{2\pi\hb}\int_{\CX}\Pi\,dm_\nu=I^0,
$$
where $I^0$ is the unity operator in the invariant subspace
$L^0\subset L$ which is generated 
by the representation of the algebra~{\rm(2.3)} 
from the fiducial state~$\CP^0$.

{\rm(c)} The mapping 
$$
\CL_\nu\to L^0\subset L,\qquad 
\psi\mapsto\int_{\CX}\frac{\psi(\oz)}{\CK_\nu(\oz|z)}\CP_z\,dm_\nu(\oz|z)
\tag 4.11
$$
defines a unitary isomorphism between Hilbert spaces.
The inverse mapping is given by the formula
$$
L\to\CL_\nu,\qquad \CP\mapsto(\CP,\CP_z).
$$

{\rm(d)} The transformation {\rm(4.11)}
determined by the coherent states~{\rm(4.10)}
intertwines the representation of the algebra~{\rm(2.3)} 
in the space~$L$ 
with the irreducible representation~{\rm(3.3)} 
in the space~$\CL_\nu$ of antiholomorphic functions
on the cylinder~$\CX$. 
\endproclaim

Now let us consider the case of a torus. 
{\it In the nonresonance version}~(3.11) 
{Theorem~{\rm 4.1} remains valid}, 
but in this case the integral in statements~(b) and~(c)
must be taken over an infinity sheet covering of the torus~$\CX$.  

Let us describe the reproducing kernel 
and the coherent states in the resonance case~(3.12). 
The space of antiholomorphic functions $\CL^N_\nu$ 
is specified by (3.16) and (3.17).
Hence we take the following orthonormal basis in~$\CL^N_\nu$:
$$
e^{(n)}(\oz)=\nu_{!}(n\hb)^{-1}
\exp\Big\{-\frac{\tau\hb n^2}{2}+n\oz\Big\}
\th \Big(\frac{N(\oz-\tau\hb n)}{i},
\frac{\tau\hb N^2}{2}\Big),
\quad n=0,1,\dots,N-1.
$$
This is the eigenbasis for the operator
$\exp\big\{\frac{2\pi i}{N\hb}\wh{t}\big\}=
\exp\big\{\frac{2\pi i}{N}\opa\big\}$. 
After the quantum Fourier transformation
$e^{(n)}\to{\widetilde e\,}^{(n)}$,
we obtain the new basis 
$$
{\widetilde e\,}^{(n)}(\oz)=
\frac{1}{\sqrt{N}}\sum^{N-1}_{n=0}e^{\frac{2\pi i}{N} mn} e^{(m)}(\oz)=
\frac{1}{\sqrt{N}}\th_{\nu}\Big(\frac{\oz}{i}+\frac{2\pi n}{N},
\frac{\tau\hb}{2}\Big),
$$
which is the eigenbasis for the operator 
$\exp\{i\wh{s}\}=\nu(\hb\opa)^{-1}\exp\{\oz-\tau\hb\opa\}$ 
from~(3.1).
Then the reproducing kernel is determined by the formula
$$
\CK^N_\nu(\oz|z)
=\sum^{N-1}_{n=0}|e^{(n)}(\oz)|^2
=\sum^{N-1}_{n=0}|{\widetilde e\,}^{(n)}(\oz)|^2
=|\nu_{!}(\hb\opa)|^{-2}\CK^N(\oz|z),
\tag 4.12
$$
where
$$
\CK^N(\oz|z)=\frac{1}{N}\sum^{N-1}_{n=0}
\Big|\th\Big(\frac{\oz}{i}+\frac{2\pi n}{N},
\frac{\tau\hb}{2}\Big)\Big|^2.
$$
This is a modification of formula (4.2) 
for the case of a resonant torus.
Note that 
the periodicity condition~(3.13)
allowed us to take the factor $\nu_{!}$ 
outside the sum symbol in~(4.12).

\proclaim{Lemma 4.2}
In the special case {\rm(3.7a)} 
\rom(i.e., for $\nu=1$ and $g=0$\rom)
the following representation for the reproducing kernel 
$\CK^N$ over the resonant torus holds\rom: 

--- if $N$ is odd then
$$
\CK^N(\oz|z)=\th \Big(\frac{\oz+z}{i},\tau\hb\Big)
\th \Big(\frac{N(\oz-z)}{2i},\frac{\tau\hb N^2}{4}\Big);
$$

--- if $N$ is even then
$$
\align
\CK^N(\oz|z)&=\th \Big(\frac{\oz+z}{i},\tau\hb\Big)
\th \Big(\frac{N(\oz-z)}{i},\tau\hb N^2\Big)\\
&\qquad 
+\th^{\#}\Big(\frac{\oz+z}{i},\tau\hb\Big)
\th^{\#}\Big(\frac{N(\oz-z)}{i},\tau\hb N^2\Big),
\endalign
$$
where 
$$
\th^{\#}(\alpha,\ve)\od\sum_{n\in\frac12+\BZ}
\exp\{-\ve n^2+in\alpha\},\qquad \ve>0.
$$
\endproclaim

Thus in the case of a resonant torus,
the reproducing kernel is determined by the product  
of two theta-functions
(a similar assertion was proposed as a hypothesis in~\cite{16}, 
but for a different space of antiholomorphic functions).  

Note that, just as in~(4.3), 
the reproducing kernel can be specified 
as a unique solution of a system of several (difference)
equations. 
Then the function $\CK^N_\nu(\oz|w)$ is a solution of the
problem 
$$
\gather
|\nu|^{-2}\Big(\hb\frac{\pa}{\pa\oz}\Big)\exp\Big\{\oz+w-
\tau\hb\Big(\frac{\pa}{\pa\oz}+\frac{\pa}{\pa w}\Big)\Big\}\CK^N_\nu=
\CK^N_\nu,
\\
\CK^N_\nu(\oz+2\pi i\,|w)=\CK^N_\nu(\oz|w), 
\qquad
\exp\Big\{N\Big(\oz-\tau\hb\frac{\pa}{\pa\oz}\Big)\Big\}\CK^N_\nu=\CK^N_\nu,
\\
\CK^N_\nu\Big(\oz+\frac{2\pi i}{N}\Big|w\Big)=
\CK^N_\nu\Big(\oz\Big|w+\frac{2\pi i}{N}\Big),
\qquad
\frac1{(2\pi)^2}\int^{2\pi}_0 d\alpha\int^{2\pi}_0 d\beta\,\,
\CK^N_\nu(i\alpha|i\beta)=1.
\endgather
$$

Next, just as in the case of a cylinder, 
we present the reproducing kernel in the form 
$$
\CK^N_\nu(\oz|z)=\sqrt{\frac{\pi}{\tau\hb}}
\exp\Big\{\frac{(\oz+z)^2}{4\tau\hb}\Big\}q^N_\nu(\oz|z)
\tag 4.13
$$
and introduce the measure on the torus $\CX$ as
$$
dm^N_\nu\od q^N_\nu(\oz|z)\, p_\nu(\oz+z)
\frac{d\oz dz}{2\tau}\,\,,
\tag 4.14
$$
where $p_{\nu}$ is the function~(4.4).

\proclaim{Lemma 4.3}
If the function $p_{\nu}$ is bounded, the norm~{\rm(3.17)} 
in the Hilbert space $\CL^N_\nu$ can be written 
in the integral form as 
$$
\align
\|\psi\|^2&=\frac1{2\pi}
\int_{\Sb 0\leq\Im z\leq 2\pi\\ 0\leq \Re z\le \tau T\endSb}
|\psi(\oz)|^2p_\nu(\oz+z)
\frac{\exp\{-(\oz+z)^2/4\hb\tau\}}{\sqrt{4\pi\hb\tau}}\,d\oz dz\\
&
=\frac{1}{2\pi\hb}\int_{\CX}\frac{|\psi|^2}{\CK^N_\nu}dm^N_\nu.
\tag 4.15
\endalign
$$
The corresponding quantum K\"ahler form 
on the resonant torus $\CX$ has the form
$$
\omega^N_{\nu}\od i\hb \opa\pa(\ln\CK^N_\nu)\,d\oz \wedge dz.
\tag 4.16
$$
The quantization condition is satisfied\rom:
$$
\frac{1}{2\pi\hb}\int_{\CX}\omega^N_{\nu}=N.
\tag 4.17
$$
In formulas {\rm (4.15)} and {\rm (4.17)} 
the operation $\int_{\CX}$ is understood, in general, 
as an integral over the $m$-sheet covering of the torus~$\CX$ 
{\rm(}see the comments on formula~{\rm (3.12))}.
\endproclaim

This is an analog of Lemma~4.1 and formula~(4.6). 
Note that the quantization condition~(4.17) 
for the K\"ahler form $\omega^N_{\nu}$ 
does not differ from condition~(3.12) 
for the form $\omc=\lim_{\hb\to 0}\omega^N_{\nu}$
since the first Chern class of tori is trivial:
$c_1(\CX)=0$. 
If the surface~$\CX$ is homeomorphic to the sphere,
then there is a distinction;
for details see~\cite{23}. 

Now we can again follow the scheme for calculating 
coherent states.
The fiducial state $\CP^0\in\CL^N_\nu$ is 
$$
\CP^0(\oz)=\th\Big(\frac{N\oz}{i},\frac{\tau\hb N^2}{2}\Big).
$$
For any $T$-periodic function~$f(t)$ the state $\CP^0$
is an eigenstate of the operator $f(\wh{t}\,)$.
Namely $f(\wh{t}\,)\CP^0=f(0)\cdot\CP^0$.
Hence it follows from (3.5) and (3.6) that
$$
\align
&\exp\Big\{-\frac{\tau\hb n^2}{2}+n\oz\Big\}
\th\Big(\frac{N(\oz-\tau\hb n)}{i},\frac{\tau\hb N^2}{2}\Big)
=e^{n(g(\wh{t})+i\wh{s}\,)}\CP^0(\oz)\\
&\qquad 
=e^{in\wh{s}}
\exp\bigg\{\frac1{\hb}\int^{\wh{t}+n\hb}_{\wh{t}}g(t)\,dt\bigg\}\CP^0
=e^{in\wh{s}}\nu_{!}(n\hb+\wh{t})\CP^0
=\nu_{!}(n\hb)e^{in\wh{s}}\CP^0.
\endalign
$$
Substituting this expression into the formula for $e^{(n)}(\oz)$ 
and then into (4.12), we obtain
$$
\align
\CK^N_\nu
&=\onu_{!}(\hb\pa)^{-1}\sum^{N-1}_{n=0}
\exp\Big\{-\frac{\tau\hb n^2}{2}+n(z+i\wh{s})\Big\}
\th\Big(\frac{N(z-\tau\hb n)}{i},\frac{\tau\hb N^2}{2}\Big)\CP^0\\
&=\onu_{!}(\hb\pa)^{-1}
\th\Big(\frac{z}{i}+\wh{s},\frac{\tau\hb}{2}\Big)\CP^0.
\endalign
$$
This formula determines coherent states in~$\CL^N_\nu$.
We see that this is identically the same formula as~(4.9).
Thus the coherent states over the resonant torus~$\CX$ 
in an abstract Hilbert space~$L$ 
are given by (4.10) if, along with~(4.8), it is required that 
$$
\wh{B}^N \CP^0=\beta\CP^0, \qquad \beta=\const.
$$
The value of~$\beta$ can be found from~(3.18):
$\beta=\CF_{!}(N\hb)^{1/2}e^{i\alpha}$. 

We have the following analog of Theorem~4.1. 

\proclaim{Theorem 4.2}
Under the assumptions of Theorem~{\rm 3.2}
the following assertions hold.

{\rm(a)} The inner product of coherent states
{\rm(4.10)} over the resonant torus~$\CX$ 
implies the reproducing kernel
$\CK^N_\nu(\oz|z)=(\CP_z,\CP_z)$.

{\rm(b)} The following partition of unity holds{\rm:}
$$
\frac1{2\pi\hb}\int_{\CX}\Pi\,dm^N_\nu=I^0,
$$
where $I^0$ is the unity operator 
in the invariant subspace $L^0\subset L$
generated by the representation of the algebra~ {\rm(2.3)}
from the fiducial state~$\CP^0$.

{\rm(c)} The mapping 
$$
\CL^N_\nu\to L^0\subset L,\qquad 
\psi\mapsto\int_{\CX}\frac{\psi(\oz)}{\CK^N_\nu(\oz|z)}\CP_z\,
dm^N_\nu(\oz|z)
\tag 4.18
$$
defined a unitary isomorphism between Hilbert spaces.
The inverse mapping is given by the formula
$$
L\to\CL^N_\nu,\qquad \CP\mapsto(\CP,\CP_z).
$$

{\rm(d)} The transformation~{\rm(4.18)} 
determined by the coherent
states~{\rm(4.10)} intertwines
the representation of the algebra~{\rm(2.3)} 
in the space~$L$ with the irreducible representation~{\rm(3.3)} 
in the space~$\CL^N_\nu$ of antiholomorphic sections over the
torus~$\CX$. 
\endproclaim

\head 5. Examples: Sklyanin algebra and algebra $su(1,1)$\endhead

\example{Example 5.1}
We consider four Hermitian generators satisfying the quadratic
relations 
$$
\alignat2
[\wh{S}_1,\wh{S}_2]&=i(\wh{S}_0\wh{S}_3+\wh{S}_3\wh{S}_0),
&\qquad [\wh{S}_0,\wh{S}_1]&=-i r^2(\wh{S}_2\wh{S}_3+\wh{S}_3\wh{S}_2),\\
[\wh{S}_2,\wh{S}_3]&=i(\wh{S}_0\wh{S}_1+\wh{S}_1\wh{S}_0),
&\qquad [\wh{S}_0,\wh{S}_2]&=i r^2(\wh{S}_3\wh{S}_1+\wh{S}_1\wh{S}_3),\\
[\wh{S}_3,\wh{S}_1]&=i(\wh{S}_0\wh{S}_2+\wh{S}_2\wh{S}_0),
&\qquad [\wh{S}_0,\wh{S}_3]&=0.
\endalignat
$$
In the paper~\cite{16} this algebra is numbered as 
the ``degenerate case~(2a).''
We assume that $r>0$ and introduce the number 
$q=\frac{1+ir}{1-ir}=e^{i\varphi}$, 
where $r=\tg\frac{\varphi}{2}$.
We also introduce the new generators 
$$
\wh{A}=\sqrt{r}\,\wh{S}_3+\frac{i}{\sqrt{r}}\wh{S}_0,\qquad
\wh{B}=\wh{S}_1-i\wh{S}_2,\qquad
\wh{C}=\wh{S}_1+i\wh{S}_2
$$
satisfying the relations 
$$
\gathered
[\wh{C},\wh{B}]=-i({\wh{A}\,}^2-\wh{A}^{*2}),
\qquad [\wh{A},\wh{A}^*]=0,\\
\wh{C}\wh{A}=q\wh{A}\wh{C},\qquad
\wh{A}\wh{B}=q\wh{B}\wh{A},\qquad
\wh{B}^*=\wh{C}.
\endgathered
\tag 5.1
$$
Here, to simplify the notation, 
we use the non-Hermitian generator 
$\wh{A}=\wh{A}_1+i\wh{A}_2$ 
instead of its real and imaginary parts.

These relations have the form~(2.3) where 
we need to set $\hb=1$ and specify the flow~$\Phi_t$
by the formula:
$\Phi_t(A_0,A)=\big(A_0
+\frac{\oq(q^{2t}-1)A^2+q(\oq^{2t}-1)\oA^2}{i(q-\oq)},\,\,\,q^tA\big)$.
Here the Casimir elements of the form (2.5) 
are given by the functions
$\varkappa_0=A_0-\frac{\oq A^2+q\oA^2}{i(q-\oq)}$
and
$\varkappa_1=A\oA$.

Now we assume that the parameters~$a_0$ and~$a$
of the surface $\CX$ (2.8) are chosen so that
$$
a_0>\varkappa_1\frac{1-\cos(\psi-\varphi)}{\sin\varphi},\qquad 
\text{where}\quad a=\sqrt{\varkappa_1}e^{i\psi/2}.
\tag 5.2
$$
Then the function 
$$
\CF(t)=a_0+\frac{\varkappa_1}{\sin\varphi}
\Big(\cos(\psi-\varphi)-\cos(\psi+(2t-1)\varphi)\Big)
$$
is strictly positive for all~$t$, 
and, obviously, the periodicity condition~(2.9) 
holds with period $T=2\pi/\varphi$. 
Hence the surface~$\CX$ is embedded in~$\BR^4$ 
as a torus.
However, the quantization condition~(3.12) (for $\hb=1$)
holds if and only if
$$
\frac{\varphi}{2\pi}\quad \text{is rational}.
\tag 5.3
$$
This implies the following condition 
on the structural constant~$r$ 
in the original commutation relations:
if $\frac1{\pi}\arctg r$ is an irrational number, 
then 
the quantization condition does not hold on any torus, 
i.e.,
there are no resonant tori;
but if $\frac1{\pi}\arctg r$ is rational, 
then any torus is a resonant torus.

Following~(3.2), 
we would like to introduce a function $\mu(t)$ so that 
$\CF(t)=|\mu(t)|^2$.
We want to avoid square roots of the form $\CF(t)^{1/2}$. 
Let us note that 
$$
\CF(t)=a_0+v_0(\varphi_t(a))-v_0(a),\qquad \text{where}\quad
v_0(A)=\frac{\oq A^2+q\oA^2}{i(q-\oq)},\quad
\varphi_t(a)=q^t a.
$$
Thus we can set $\mu(t)=M(\varphi_t(a))$, where
$M(A)\overline{M(A)}=a_0-v_0(a)+v_0(A)$.
Let us seek the function~$M$ in the form
$M(A)=\zeta A-\overline{\xi}\oA$, 
where $|A|=|a|\equiv\sqrt{\varkappa_1}$.
Then we obtain the following system for the coefficients~$\zeta$
and~$\xi$ 
(see the notation in~(5.2)): 
$$
\zeta\xi=\frac{e^{-i\varphi}}{2\sin\varphi},\qquad
|\zeta|^2+|\xi|^2=\frac{a_0}{\varkappa_1}
+\frac{\cos(\psi-\varphi)}{\sin\varphi}.
\tag 5.4
$$
This system is easily solved.
So we choose 
$$
\mu(t)=\zeta a\,e^{i\varphi t}-\overline{\xi a}\,e^{-i\varphi t},
\tag 5.5
$$
where $\zeta$ and $\xi$ are subject to~(5.4).

Let us take the simplest factor~$\CB(t)$ in~(3.4): $\CB=\mu$. 
Then $\nu=1$, $g=0$, and the complex structure on~$\CX$ 
is determined by~(3.7a): 
$\oz=\tau t+is$, where $\tau>0$.

\smallskip

(1) {\tt Nonresonance version}: (5.3) does not hold.
Here it is necessary to consider 
a cylinder infinitely wound on the torus~$\CX$. 
In this case, by Lemma~4.1, 
the Hilbert space of $2\pi i$-periodic antiholomorphic functions
over the covering of the torus~$\CX$
is endowed with the norm
$$
\|\psi\|^2=\frac1{2\pi}\int_{0\leq\Im z\leq 2\pi}
|\psi(\oz)|^2\frac{\exp\{-(\oz+z)^2/4\tau\}}{\sqrt{4\pi\tau}}
d\oz\,dz.
\tag 5.6
$$
By~(4.2), 
the reproducing kernel of this space is determined 
by the theta-function:
$\CK(\oz|z)=\th\big(\frac{\oz+z}{i},\tau\big)$.
The quantum K\"ahler form and the reproducing 
measure are given by formulas~(4.7) and, as is easily seen,
are well defined only on the infinite sheet covering 
of the torus (on the cylinder).

The irreducible representation of relations~(5.1) 
in the Hilbert space~(5.6) is specified by operators of the
form~(3.3): 
$$
\alignedat2
&\wh{A}=a\exp\{i\varphi\opa\},&\qquad 
&\wh{A}^*=\overline{a}\exp\{-i\varphi\opa\},\\
&\wh{B}=(\zeta\wh{A}-\overline{\xi}\wh{A}^*)\exp\{\oz-\tau\opa\},
&\qquad 
&\wh{C}=\exp\{\tau\opa-\oz\}(\overline{\zeta}\wh{A}^*-\xi\wh{A}).
\endalignedat
\tag 5.7
$$

(2) {\tt Resonance version}: condition (5.3) is satisfied, 
i.e., $\varphi=2\pi m/N$, where $m$ and $N$ are coprime integers.
For the period we choose $T=N$ (the minimal period is equal
to~$N/m$).   
Then over the $m$-multiple covering of the torus 
we construct the $N$-dimensional Hilbert space 
of functions satisfying the quasiperiodicity condition (3.14) 
(where $\hb=1$).
This Hilbert space is endowed with the norm 
$$
\|\psi\|^2=\frac1{2\pi}\int_{\Sb 0\leq\Im z\leq 2\pi\\
0\leq\Re z\leq \tau T\endSb}|\psi(\oz)|^2
\frac{\exp\{(\oz+z)^2/4\tau\}}{\sqrt{4\pi\tau}}
d\oz\,dz.
$$
The reproducing kernel of this space 
is given in Lemma~4.2 (where $\hb=1$).
Here the form $\omega^N_{\nu}$ 
and the measure $dm^N_{\nu}$ are determined by 
(4.13), (4.14), and~(4.16) 
with $\nu=1$ as geometric objects on the torus,
more precisely, on its $m$-multiple covering.

Since we have $\nu=1$ in this case, 
the first normalization condition~(3.13) 
holds automatically.
The second condition~(3.13) can be ensured as follows:
it is necessary to replace 
the originally chosen solution $\zeta$, $\xi$ of system (5.4) 
by another solution 
$\widetilde{\zeta}$, $\widetilde{\xi}$ 
according to the formulas
$$
\widetilde{\zeta}=\zeta\exp\Big\{i\frac{\alpha-\delta}{N}\Big\},\quad
\widetilde{\xi}=\xi\exp\Big\{i\frac{\delta-\alpha}{N}\Big\},
$$
where $\alpha$ is the parameter from~(3.13), 
$\delta=\sum^{N}_{n=1}\arg\mu(n\hb)$, 
and the function~ $\mu$ is defined in~(5.5) 
by using the solution $\zeta$, $\xi$.
The operators of the irreducible representation of the
algebra~(5.1) are determined by the same formulas~(5.7)
with $\zeta$ and $\xi$ replaced by $\widetilde{\zeta}$ 
and $\widetilde{\xi}$.

Note that, in the resonance version, 
the algebra~(5.1), 
in addition to two ``classical'' Casimir elements,
also possesses two ``nonclassical'' elements 
$\wh{B}^N$, $\wh{A}^N$ and their adjoints 
(which are scalars in the operator irreducible 
representation). 

\example{Remark 5.1}
In~\cite{16} the nonresonance version was not studied.
It should be noted that in this version 
infinite-dimensional representations are assigned to 
compact symplectic leaves of the corresponding 
Poisson algebra.
In the resonance version, our representations~(5.7) 
defined on antiholomorphic functions 
were also not studied in~\cite{16} 
(for this case the representations in~\cite{16}
are constructed in the space of functions of a circle).
\endexample
\endexample

\example{Example 5.2}
Now we consider the Lie algebra $su(1,1)$. 
Its three Hermitian generators satisfy 
the commutation relations
$$
[\wh{S}_1,\wh{S}_2]=i\hb\wh{S}_3,\qquad 
[\wh{S}_2,\wh{S}_3]=-i\hb\wh{S}_1,\qquad 
[\wh{S}_3,\wh{S}_1]=-i\hb\wh{S}_2.
$$
We denote $\wh{B}=\wh{S}_1-i\wh{S}_2$, 
$\wh{C}=\wh{S}_1+i\wh{S}_2$, and $\wh{A}=\wh{S}_3$. 
Then the relations become
$$
\wh{C}\wh{B}=\wh{B}\wh{C}+2\hb\wh{A},\qquad 
\wh{C}\wh{A}=(\wh{A}+\hb)\wh{C},\qquad 
\wh{B}^*=\wh{C},\qquad 
\wh{A}^*=\wh{A}.
\tag 5.8
$$
This is a special case of relations (2.3)
where the flow $\Phi_t:\BR^2\to\BR^2$ 
has the form
$\Phi_t(A_0,A)=\big(t^2+t(2A-\hb)+A_0,\,\,A+t\big)$.
Assume that the parameters~$a_0$ and~$a$ are chosen so that 
$a_0-(a-{\hb}/{2})^2\od\lambda^2>0$.
Then the function (3.2) is positive:
$\CF(t)=t^2+t(2a-\hb)+a_0>0$ for all $t\in\BR$.
Thus the surface $\CX$ is diffeomorphic to a cylinder
embedded in $\BR^3$ as a one-sheet hyperboloid 
$$
\CX=\big\{BC-(A-{\hb}/{2})^2=\lambda^2\big\}.
$$
We choose the function $\mu(t)$ from (3.2) as
$\mu(t)=t+a-{\hb}/{2}-i\lambda$.
Now we consider two versions of choosing the factor~$\CB$
in~(3.4). 
Namely, we choose either $\CB(t)\equiv\mu(t)$ or
$\CB(t)\equiv\CF(t)$. 

\smallskip

{\tt Version I}. Let us choose $\CB=\mu$. 
Then $\nu=1$, and the Hilbert space of antiholomorphic functions
on the cylinder is determined by the norm 
(4.5) (with $p_\nu\equiv 1$). 
The irreducible representation of the Lie algebra $su(1,1)$ 
(i.e.,  the representation of relations~(5.8)) 
is given by the following operators acting in this Hilbert space:
$$
\wh{A}=a+\hb\opa,\quad
\wh{B}=(\wh{A}-{\hb}/{2}-i\lambda)\,
e^{\oz-\tau\hb\opa},\quad
\wh{C}=e^{\tau\hb\opa-\oz}\,
(\wh{A}-{\hb}/{2}+i\lambda).
\tag 5.9
$$

{\tt Version II}. Let us choose $\CB=\CF$.
Then $\nu={\overline{\mu}\,}^{-1}$, 
and the complex structure on~$\CX$ is determined by~(3.7)
with 
$$
g(t)=-\frac{\hb}{a+t-\frac{\hb}2+i\lambda}-
\hb\frac{d}{dt}\ln\Gamma\Big(\frac{a+t}{\hb}-\frac12
+\frac{i\lambda}{\hb}\Big)-\ln\hb,
$$
where $\Gamma$ is the standard gamma-function.
Hence we have
$$
\nu_{!}(t)=\frac{\Gamma(\frac{a}{\hb}+\frac12+\frac{i\lambda}{\hb})}
{\Gamma(\frac{a+t}{\hb}+\frac12+\frac{i\lambda}{\hb})}
\exp\Big\{-\frac{t}{\hb}\ln\hb\Big\},
$$
which implies the following formula for the function (4.4):
$$
\multline
p_\nu(x)=
\frac{\Big|\Gamma\Big(\frac{a}{\hb}+\frac12+\frac{i\lambda}{\hb}\Big)\Big|^2}
{\sqrt{\pi\hb\tau}}\\
\times\int^{\infty}_{-\infty}
\frac{\exp\{-\frac{t^2}{\tau\hb}\}\exp\Big\{-\frac{x+2it}{\tau\hb}\ln\hb\Big\}}
{\Gamma(\frac{a}{\hb}+\frac12+\frac{x}{2\tau\hb}
+\frac{it}{\tau\hb}+\frac{i\lambda}{\hb})\,
\Gamma(\frac{a}{\hb}+\frac12+\frac{x}{2\tau\hb}
+\frac{it}{\tau\hb}-\frac{i\lambda}{\hb})}\,dt.
\endmultline
\tag 5.10
$$
Thus, in this version, the Hilbert space~$\CL_\nu$
of antiholomorphic functions on the cylinder is given by the
norm~(4.5), where $p_\nu$ is defined in~(5.10).
The reproducing kernel of this space is given 
by the modified theta-function~(4.2):
$$
\CK_\nu(\oz|z)
=\frac{|\Gamma(\frac{a}{\hb}+\frac12+\frac{i\lambda}{\hb}+\pa)|^2}
{|\Gamma(\frac{a}{\hb}+\frac12+\frac{i\lambda}{\hb})|^2}
\th\Big(\frac{\oz+z+2\ln\hb}{i},\tau\hb\Big).
$$
The irreducible representations of the Lie algebra $su(1,1)$ 
in the Hilbert space $\CL_\nu$ has the form (3.3):
$$
\wh{A}=a+\hb\opa,\qquad
\wh{B}=\big[(\wh{A}-{\hb}/{2})^2+\lambda^2\big]
e^{\oz-\tau\hb\opa},\qquad
\wh{C}=e^{\tau\hb\opa-\oz}.
\tag 5.11
$$

\example{Remark 5.2}
Apparently, 
the representations (5.9) and (5.11) have not been studied 
in the standard representation theory of Lie algebras. 
They are associated with the complex structure
(the complex polarization), 
which is not invariant under the co-adjoint action of~$su(1,1)$
on the symplectic leaf~$\CX$. 
The usual way is to study representations corresponding to 
the invariant real polarization of~$\CX$
(the fibration by circles $A=\const$).
\endexample

\end